\documentclass[12pt]{amsart}
\usepackage{a4}
\usepackage{amsfonts, amssymb}

 \address{Instituto de Ciencias Matem\'{a}ticas (CSIC-UAM-UCM-UC3M) and Departamento de Matem\'aticas \\ Universidad Aut\'onoma de Madrid \\ 28049 Madrid, Spain}

  \email{franciscojavier.cilleruelo@uam.es}

\newtheorem{thm}{Theorem}
\newtheorem{cor}{Corollary}

\newtheorem{lem}{Lemma}

\newtheorem{defn}{Definition}

\numberwithin{equation}{section} \numberwithin{thm}{section}
\numberwithin{lem}{section} \numberwithin{problem}{section}
\numberwithin{cor}{section}

\newcommand{\p}{\mathbf p}

\newcommand{\Z}{\mathbb{Z}}

\begin{document}

\title{The least common multiple of a  quadratic sequence}
\author{Javier Cilleruelo}
\thanks{This work was supported by Grant MTM 2008-03880 of MYCIT (Spain)}

\begin{abstract}
For any irreducible quadratic polynomial $f(x)$ in $\Z[x]$ we obtain
the estimate $\log \ \text{l.c.m. } \{f(1),\dots ,f(n)\}= n\log
n+Bn+o(n)$ where $B$ is a constant depending on $f$.
\end{abstract}

\maketitle \small
\section{Introduction}It is well known that $\log \text{l.c.m.}\{1,\dots , n\}\sim n $. Indeed, this asymptotic estimate is equivalent to the prime number theorem.
The analogous for arithmetic progressions is also known
\cite{Bateman} and it is a consequence of the prime number theorem
for arithmetic progressions:
\begin{equation}\label{Ba}\log \text{l.c.m.}\{a+b,\dots ,
an+b\}\sim n\frac{q}{\phi (q)}\sum_{\substack{1\le k\le q\\
(k,q)=1}}\frac 1k, \end{equation} where $q=a/(a,b)$.

We address here the problem of estimating $\log \text{l.c.m.}
\{f(1),\dots ,f(n)\}$ when $f$ is an irreducible quadratic
polynomial in $\Z[x]$. The same problem for reducible quadratic
polynomials  is  easier and we study it in section \S4.
%

\smallskip

\begin{thm}\label{principal} For any irreducible quadratic polynomial
$f(x)=ax^2+bx+c$ in $\Z[x]$ we have
$$\log  \text{l.c.m. }\{f(1),\dots ,f(n)\}= n\log n +Bn+o(n)$$
where  $B=B_f$ is defined by the formula
\begin{align}\label{B}
B_f&=\gamma-1-2\log 2-\sum_{p}\frac{(d/p)\log p}{p-1}+\frac
1{\phi(q)}\sum_{\substack{1\le r\le q\\(r,q)=1}}\log\Bigl (1+\frac
rq\Bigr )\\&+\log a+\sum_{p\mid 2aD}\log p\Bigl(\frac{1+(d/p)}{p-1}-
\sum_{k\ge 1}\frac{s(f,p^k)}{p^k}\Bigr ).\notag
\end{align}
In this formula $\gamma$ is the Euler constant, $D=b^2-4ac$, $d$ is
the fundamental discriminant, $(d/p)$ is the Kronecker symbol,
$q=a/(a,b)$ and $s(f,p^k)$ is the number of solutions of $f(x)\equiv
0\pmod{p^k}$, which can be calculated easily using lemma \ref{H}.
\end{thm}
In section \S 3 we give an alternative expression for the constant
$B_f$, which is more convenient for numerical computations. As an
example we will see that for the simplest case, $f(x)=x^2+1$, the
constant $B_f$ in theorem \ref{principal} can be written as
\begin{align*}\label{simplest}B_f&=\gamma -1-\frac{\log
2}2-\sum_p\frac{(-4/p)\log p}{p-1}\\&=\gamma -1-\frac{\log
2}2+\sum_{k=1}^{\infty}\frac{\zeta'(2^k)}{\zeta(2^k)}+\sum_{k=0}^{\infty}\frac{L'(2^k,\chi_{-4})}{L(2^k,\chi_{-4})}-\sum_{k=1}^{\infty}\frac{\log
2}{2^{2^k}-1}
\\&=-0.066275634213060706383563177025...\end{align*}

\smallskip

It would be interesting to extend our estimates to irreducible
polynomials of higher degree, but we have found a serious
obstruction in our argument. Some heuristic arguments and
computations allow us to conjecture that the asymptotic estimate
$$\log
 \text{l.c.m. }\{f(1),\dots ,f(n)\}\sim (\text{deg}(f)-1)n\log n$$ holds
for any irreducible polynomial $f$ in $\Z[x]$ of degree $\ge 2$.

An important ingredient in the proof of theorem \ref{principal} is a
deep result about the distribution of the solutions of the quadratic
congruences $f(x)\equiv 0\pmod p$ when $p$ runs over all the primes.
It was proved by Duke, Friedlander and Iwaniec \cite{DFI} (for $D
<0$), and by Toth (for $D >0$). Actually we will need a more general
statement of this result, due to Toth.
\begin{thm}[Toth \cite{T}]\label{I}
For any irreducible quadratic polynomial $f$ in $\Z[x]$, the
sequence $$\{\nu/p,\ 0\le \nu <p\le x,\ p\in S, \ f(\nu)\equiv
0\pmod p\}$$ is well distributed in $[0,1)$ as $x$ tends to infinity
for any arithmetic progression $S$ containing infinitely many primes
$p$ for which the congruence $f(x)\equiv 0\pmod p$ has solutions.
\end{thm}

\

\textbf{Acknowledgment.} We thank to Arpad Toth for clarifying  the
statement of theorem 1.2 in \cite{T}, to Guoyou Qian for detecting a
mistake in Lemma \ref{H}, to Adolfo Quir\'{o}s for conversations on some
algebraic aspects of the problem, to Enrique Gonz\'{a}lez Jim\'{e}nez for
the calculations of some constants and to Fernando Chamizo for some
suggestions and a carefully reading of the paper.

\section{Proof of theorem \ref{principal}}
\subsection{Preliminaries}
For $f(x)=ax^2+bx+c$ we define $D=b^2-4ac$ and
$$L_n(f)=\text{l.c.m.} \{f(1),\dots ,f(n)\}.$$
 Since $L_n(f)=L_n(-f)$ we can assume that $a>0$. Also we can assume that $f(x)$ is positive and increasing for $x\ge 1$. If it is not the case,
 we consider a polynomial $f_k(x)=f(k+x)$ for a  $k$ such that $f_k(x)$ is positive and increasing for $x\ge 1$. Then we observe that
 $L_n(f)=L_n(f_k)+O_k(\log n)$ and that the error term is negligible for the statement of theorem \ref{principal}.

We define the numbers $\beta_p(n)$ by the formula
\begin{equation}\label{Ln}L_n(f)=\prod_pp^{\beta_p(n)}\end{equation}where the
product runs over all the primes $p$. The primes involved in this
product are those for which the congruence $f(x)\equiv 0\pmod p$ has
some solution. Except for some special primes (those such that
$p\mid 2aD$) the congruence $f(x)\equiv 0\pmod {p}$ has $0$ or $2$
solutions. We will discus it in detail in lemma \ref{H}.

We denote by $\mathcal P_f$ the set of the non special primes for
which the congruence $f(x)\equiv 0\pmod p$ has exactly two
solutions. More concretely
$$\mathcal P_f=\{p:\ p\nmid 2aD,\ (D/p)=1\}$$ where $(D/p)$ is the Kronecker
symbol. This symbol is just the Legendre symbol when $p$ is an odd
prime.

The quadratic reciprocity law shows that the set $\mathcal P_f$ is
the set of the primes lying in exactly $\varphi(4D)/2$ of the
$\varphi(4D)$ arithmetic progressions modulo $4D$, coprimes with
$4D$. As a consequence of the prime number theorem for arithmetic
progressions we have
$$\# \{p\le x:\ p\in \mathcal P_f\}\sim \frac x{2\log x}$$ or equivalently,
$$\sum_{\substack{0\le \nu<p\le x\\ f(\nu)\equiv 0\pmod p}}1\sim
\frac x{\log x}.$$

Let $C=2a+b$. We classify the primes involved in \eqref{Ln}  in
\begin{itemize}
    \item \emph{Special} primes: those such that $p\mid 2aD.$

    \smallskip

\item \emph{$p\in \mathcal P_f:$}
 $\begin{cases} \text{\emph{Small} primes}:
 \ p<n^{2/3}.\\
\text{\emph{Medium} primes:} \ n^{2/3}\le p< Cn\ : \begin{cases}bad
\text{ primes: }p^2\mid f(i) \text{ for some } i\le n. \\good \text{
primes: } p^2\nmid f(i) \text{ for any } i\le n.  \end{cases}\\
\text{\emph{Large} primes:} \ Cn\le p\le f(n).
\end{cases}$
    \end{itemize}

 We will use different strategies to deal with these primes.

\subsection{Large primes}
To deal with the large primes we consider
 $P_n(f)$ and the numbers $\alpha_p(n)$ defined by the
expression
\begin{eqnarray}
P_n(f)&=&\prod_{i=1}^nf(i)=\prod_p p^{\alpha_p(n)}.
\end{eqnarray}
Next lemma allow us to avoid the large primes.
\begin{lem}
If $p\ge 2an+b$ then $\alpha_p(n)=\beta_p(n)$.
\end{lem}
\begin{proof}If $\beta_p=0$ then $\alpha_p(n)=0$.
If $\alpha_p(n)>\beta_p(n)\ge 1$  there exist $i<j\le n$ such that
$p\mid f(i)$ and $p\mid f(j)$. It implies that $p\mid
f(j)-f(i)=(j-i)(a(j+i)+b)$. Thus $p\mid (j-i)$ or $p\mid a(j+i)+b$,
which is not possible because $p\ge 2an+b$.
\end{proof}
Since $C=2a+b$ we can write
\begin{equation}\label{L}
\log L_n(f)=\log P_n(f)+\sum_{p<Cn}(\beta_p(n)-\alpha_p(n))\log p.
\end{equation}

Indeed we can take $C$ to be any constant  greater than $2a+b$. As
we will see, the final estimate of $\log L_n(f)$ will not depend on
$C$.

The estimate of $\log P_n(f)$ is easy: \begin{eqnarray}\log
P_n(f)&=&\log \prod_{k=1}^nf(k)=\log \prod_{k=1}^nak^2\Bigl (1+\frac
b{ka}+\frac{c}{k^2a}\Bigr )\\ &=&n\log a+\log
(n!)^2+\sum_{k=1}^n\log \Bigl (1+\frac b{ka}+\frac{c}{k^2a}\Bigr )\nonumber \\
&=&2n\log n+n(\log a-2)+O(\log n)\nonumber \end{eqnarray} and we
obtain
 \begin{equation}\label{fL}
\log L_n(f)=2n\log n+n(\log
a-2)+\sum_{p<Cn}(\beta_p(n)-\alpha_p(n))\log p+O(\log n).
\end{equation}
\subsection{The number of solutions of $f(x)\equiv 0\pmod{p^k}$ and the special primes} The number of solutions of the congruence $f(x)\equiv 0\pmod {p^k}$
will play an important role in the proof of theorem \ref{principal}.
We write $s(f,p^k)$ to denote this quantity.

Lemma belove resumes all the casuistic for $s(f,p^k)$. We observe
that except for a finite number of primes, those dividing $2aD$, we
have that $s(f;p^k)=2$ or $0$ according with $(D/p)=1$ or $-1$.
\begin{lem}\label{H}Let $f(x)=ax^2+bx+c$ be an irreducible polynomial and $D=b^2-4ac$.
\begin{enumerate}
        \item If $p\nmid 2a,\ D=p^lD_p$ and  $(D_p,p)=1$ then $$s(f,p^k)=\begin{cases}p^{\lfloor k/2 \rfloor},\ &k\le l\\ 0,\ &k> l,\ l \text{ odd or
    }(D_p/p)=-1\\2p^{l/2},\ &k> l,\ l\text{ even }\ (D_p/p)=1.
   \end{cases}$$
    \item If $p\mid a,\ p\ne 2$ then $s(f,p^k)=\begin{cases}0,\text{ if }  p\mid b\\1,\text{ if } p\nmid
    b.
    \end{cases}$
    \item If $b$ is odd then, for all $k\ge 2$, $s(f,2^k)=s(f,2)=\begin{cases}1 \text{ if } a \text{ is even}\\
    0 \text{ if } a \text{ is odd and } c \text{ is odd}\\2 \text{ if } a \text{ is odd and } c \text{ is even}.\end{cases}$
    \item If $b$ is even and $a$ is even then $s(f,2^k)=0$ for any
    $k\ge 1$.

    \item If $b$ is even and $a$ is odd, let $D=4^lD',\ D'\not \equiv 0 \pmod 4$.
\begin{enumerate}
        \item If $k\le 2l-1$, $\quad s(f;2^k)=2^{\lfloor k/2\rfloor }$
        \item  If $k=2l$, $\ \ \qquad s(f;2^k)=\begin{cases}2^{l},\ &D'\equiv 1\pmod 4  \\0,\ &D'\not\equiv 1\pmod 4. \end{cases}$
        \item  If $k\ge 2l+1$, $\quad s(f;2^k)=\begin{cases}2^{l+1},\ &D'\equiv 1\pmod 8  \\0,\ &D'\not \equiv 1\pmod 8. \end{cases}$
    \end{enumerate}

\end{enumerate}
\end{lem}
\begin{proof}
The proof is a consequence of elementary manipulations and Hensel's
lemma.
\end{proof}
\begin{cor}\label{co}If $p\nmid 2aD$ then $s(f,p^k)=1+(D/p)$.
\end{cor}
\begin{proof}
In this case, $l=0$ and $D_p=D$ in lemma \ref{H}. Thus
$s(f,p^k)=0=1+(D/p)$ if $(D/p)=-1$ and $s(f,p^k)=2=1+(D/p)$ if
$(D/p)=1$.
\end{proof}
 \begin{lem}\label{alpha}
\begin{equation}\label{alpha6}
\alpha_p(n)=n\sum_{k\ge 1}\frac{s(f,p^k)}{p^k} +O\left (\frac{\log
n}{\log p}\right ).
\end{equation}
where $s(f;p^k)$ denotes the number of solutions of $f(x)\equiv
0\pmod{p^k},\ 0\le x<p^k$.
 \end{lem}
 \begin{proof}
 We observe that the maximum exponent $\alpha_{p,i}$ such that $p^{\alpha_{p,i}}\mid f(i)$
 can be written as $\alpha_{p,i}=\sum_{k\ge 1,\ p^k\mid f(i)}1.$ Thus
 \begin{equation}\label{alpha2}
\alpha_p(n)=\sum_{i\le n}\alpha_{p,i}=\sum_{i\le n}
\sum_{\substack{k\ge 1\\ p^k\mid f(i)}}1=\sum_{k\ge 1}
\sum_{\substack{i\le n\\ p^k\mid f(i)}}1.
 \end{equation}
The trivial estimate $s(f;p^k)\left [\frac n{p^k}\right ]\le
\sum_{i\le n,\ p^k\mid f(i)}1\le s(f;p^k)\left (\left [\frac
n{p^k}\right ]+1\right )$ gives
\begin{equation}\label{s}\sum_{\substack{i\le n\\ p^k\mid
f(i)}}1=n\frac{s(f;p^k)}{p^k}+O(s(f;p^k)).\end{equation} Putting
(\ref{s}) in (\ref{alpha2}) and observing that $k\le \log f(n)/\log
p$ and that $s(f,p^k)\ll 1$, we get
\begin{align*}
\alpha_p(n)=n\sum_{k\ge 1}\frac{s(f,p^k)}{p^k}+O\Bigl ( \frac{\log
n}{\log p}\Bigr ).
\end{align*}
\end{proof}
We observe that $\beta_p(n)=\max_{i\le n}\alpha_{p,i}$, so
\begin{equation}\label{beta}\beta_p(n)\ll \log n/\log p.
\end{equation}
Now we put \eqref{beta} and \eqref{alpha6} in \eqref{fL} for the
special primes  obtaining
\begin{align}\label{Ln2}
\log L_n(f)=2n\log n&+n\left (\log a-2-\sum_{p\mid 2aD}\sum_{k\ge
1}\frac{s(f,p^k)\log p}{p^k}\right )\\&+\sum_{p<Cn,\ p\nmid
2aD}(\beta_p(n)-\alpha_p(n))\log p+O(\log n).\notag
\end{align}

\subsection{Small primes}
Lemma \ref{alpha} has an easier formulation for small primes.
\begin{lem}\label{D}For any $p\nmid 2aD$ we have
\begin{equation}\label{aa}
\alpha_p(n)=n\frac{1+(D/p)}{p-1}+O\Bigl (\frac{\log n}{\log p}\Bigr
).
\end{equation}
\end{lem}
\begin{proof}
It is a consequence of lemma \ref{alpha} and corollary \ref{co}.
\end{proof}

By substituting \eqref{aa} and \eqref{beta} in \eqref{Ln2} we obtain
\begin{align}\label{Ln3}
\log L_n(f)&=2n\log n+n\left (\log a-2-\sum_{p\mid 2aD}\sum_{k\ge
1}\frac{s(f,p^k)\log p}{p^k}\right ) \\&-\sum_{\substack{p<
n^{2/3}\\p\nmid 2aD}}\frac{(1+(D/p))\log
p}{p-1}+\sum_{\substack{n^{2/3}\le p<Cn\\ p\in \mathcal
P_f}}(\beta_p(n)-\alpha_p(n))\log p+O(n^{2/3}).\notag
\end{align}

\subsection{Medium primes}
Medium primes also can be classified in bad and good primes. Bad
primes are those $p$ such that $p^2\mid f(i)$ for some $i\le n$.
Good primes are those are not bad primes.

 As
we have seen in the previous section, for any prime $p\in \mathcal
P_f$ the congruence $f(x)\equiv 0\pmod p$ has exactly two solutions,
say $0\le \nu_{p,1},\nu_{p,2}<p$.

If $p$ is a good prime we have that $\alpha_p(n)$ is just the number
of integers $i\le n$ such that $p\mid f(i)$. These integers are all
of the form
\begin{align} &\nu_{p,1}+kp,\quad 0\le k\le \left [\frac{n-\nu_{p,1}}p\right ]\\
&\nu_{p,2}+kp,\quad 0\le k\le \left [\frac{n-\nu_{p,2}}p\right
].\end{align}

Also it is clear that if $p$ is a good prime then $\beta_p(n)\le 1$.
These observations motive the following definition:
\begin{defn}
 For any $p\in
\mathcal P_f$ we define
\begin{eqnarray}\alpha_p^*(n)&=&\left
[\frac{n-\nu_{p,1}}p\right ]+\left [\frac{n-\nu_{p,2}}p\right ]+2\\
\beta_p^*(n)&=&\begin{cases} 1,\text{ if }\beta_p(n)\ge 1\\ 0,
\text{ otherwise}.\end{cases}\end{eqnarray}
\end{defn}
\begin{lem}\label{b}For any $p\in \mathcal P_f$ we have
\begin{enumerate}
       \item[i)]$\qquad \alpha_p(n)-\alpha_p^*(n)=\frac{2n}{p(p-1)}+O(\log n/\log p)$
    \item[ii)] $\qquad \alpha_p(n)=\alpha_p^*(n)$ and
    $\beta_p(n)=\beta_p^*(n)$ if $p^2\nmid f(i)$ for any $i\le n$.
\end{enumerate}
\end{lem}
\begin{proof}
i) Lemma \ref{D} implies that $\alpha_p(n)=\frac{2n}{p-1}+O(\log
n/\log p)$ when $p\in \mathcal P_f$. On the other hand we have that
$\alpha_p^*(n)=\frac {2n}p+O(1)$. Thus,
$\alpha_p(n)-\alpha_p^*(n)=\frac{2n}{p-1}+O(\log n/\log p)-\frac
{2n}p+O(1)=\frac{2n}{p(p-1)}+O(\log n/\log p)$.

ii) The first assertion has been explained at the beginning of the
subsection. For the second, if $p\nmid f(i)$ for any $i\le n$ then
$\beta_p(n)=\beta_p^*(n)=0$. And if $p\mid f(i)$ for some $i\le n$
we have that $\beta_p^*(n)=\beta_p(n)=1$ since $p^2\nmid f(i)$.
\end{proof}

Now we split the last sum in \eqref{Ln3} in
\begin{align}\label{cc}
\sum_{\substack{n^{2/3}\le p<Cn\\ p\in \mathcal
P_f}}(\beta_p(n)-\alpha_p(n))\log p&=\sum_{\substack{n^{2/3}\le
p<Cn\\ p\in \mathcal
P_f}}(\beta_p(n)-\beta_p^*(n)-\alpha_p(n)+\alpha_p^*(n))\log p\\
&+\sum_{\substack{ p<Cn\\ p\in \mathcal
P_f}}\beta_p^*(n)\log p -\sum_{\substack{n^{2/3}\le p<Cn\notag \\
p\in \mathcal P_f}}\alpha_p^*(n)\log
p+O(n^{2/3})\\&=S_1(n)+S_2(n)-S_3(n)+O(n^{2/3}).\notag
\end{align}
To estimate $S_1(n)$ we observe that lemma \ref{b} ii) implies that
$\beta_p(n)-\beta_p^*(n)-\alpha_p(n)+\alpha_p^*(n)=0$ for any good
prime $p$. On the other hand, lemma \ref{b} i) and \eqref{beta}
implies that $|\beta_p(n)-\beta_p^*(n)-\alpha_p(n)+\alpha_p^*(n)|\ll
\log n/\log p$. Thus,
\begin{equation}
|S_1(n)|\ll \log n\ |\{p:\ n^{2/3}<p<Cn,\ p \text{ bad}\}|.
\end{equation}
\begin{lem}
The number of bad primes $p\nmid D,\ Q\le p<2Q$ is $\ll n^2/Q^2$.
\end{lem}
\begin{proof}
Let $P_r$ the set of the primes $p$ such that $f(i)=ai^2+bi+c=rp^2$
for some $i\le n$. It implies that $(2ai+b)^2-4arp^2=D$ and then,
$|\frac{2ai+b}p-2\sqrt{ ra}|\ll \frac 1{p^2}\ll \frac 1{Q^2}$. We
observe that all the fractions $\frac{2ai+b}p,\ 1\le i\le n,\ Q\le
p<2Q$ are all distinct. Otherwise $(2ai+b)p'=(2ai'+b)p$ and then
$p\mid 2ai+b$. But it would imply that $p\mid D$, which is not
possible. On the other hand $\left |\frac{2ai+b}p-\frac
{2ai'+b}{p'}\right |\ge \frac 1{pp'}\gg \frac 1{Q^2}$. Thus, the
numbers of primes $p\in P_r$ lying in $[Q,2Q]$ is $ \ll 1$. We
finish the proof observing that $r\le f(n)/Q^2\ll n^2/Q^2$.
\end{proof}

Now, if we split the interval $[n^{2/3},Cn]$ in dyadic intervals and
apply lemma above to each interval to obtain  $|S_1(n)|\ll
n^{2/3}\log n.$

To estimate $S_3(n)=\sum_{n^{2/3}<p<Cn,\ p\in \mathcal
P_f}\alpha_p^*(n)$ we start by writing
\begin{align*}\alpha_p^*(n)&=\left [\frac{n-\nu_{p,1}}p\right ]+\left
[\frac{n-\nu_{p,2}}p\right ]+2\\&=\frac{2n}p+\Bigl (\frac
12-\frac{\nu_{p,1}}p\Bigr )+\Bigl (\frac 12-\frac{\nu_{p,2}}p\Bigr
)+ \frac 12-\left\{\frac{n-\nu_{p,1}}p\right \}+\frac
12-\left\{\frac{n-\nu_{p,2}}p\right \}.\end{align*} Thus
\begin{align}
S_3(n)&=n\sum_{n^{2/3}<p<Cn}\frac{(1+(D/p))\log p}p\\
&+\sum_{\substack{n^{2/3}<p<Cn\\ 0\le \nu <p\\f(\nu)\equiv 0\pmod
p}}\Bigl (\frac 12-\frac{\nu}p\Bigr )\log
p+\sum_{\substack{n^{2/3}<p<Cn\\ 0\le \nu <p\\f(\nu)\equiv 0\pmod
p}}\Bigl (\frac 12-\Bigl \{\frac{n-\nu}p\Bigr \}\Bigr )\log p\\
&=n\sum_{n^{2/3}<p<Cn}\frac{(1+(D/p))\log p}{p-1}+O(n^{2/3})\\
&+\sum_{\substack{0\le \nu<p<Cn\\f(\nu)\equiv 0\pmod p}}\Bigl (\frac
12-\frac{\nu}p\Bigr )\log p+\sum_{\substack{0\le \nu
<p<Cn\\f(\nu)\equiv 0\pmod p}}\Bigl (\frac 12-\Bigl
\{\frac{n-\nu}p\Bigr \}\Bigr )\log p
\end{align}

Putting this in \eqref{cc} and then in \eqref{Ln3} we obtain
\begin{align}\label{Ln4}
\log L_n(f)&=2n\log n+n\left (\log a-2-\sum_{p\mid 2aD}\sum_{k\ge
1}\frac{s(f,p^k)\log p}{p^k}\right )\\& -\sum_{\substack{p<
Cn\\p\nmid 2aD}}\frac{(1+(D/p))\log
p}{p-1}+S_2(n)-T_1(n)-T_2(n)+O(n^{2/3}\log n)\notag
\end{align}
where
\begin{align}
S_2(n)&=\sum_{\substack{ p<Cn\\ p\in \mathcal
P_f}}\beta_p^*(n)\log p\label{S2}\\
T_1(n)&=\sum_{\substack{0\le \nu<p<Cn\\f(\nu)\equiv 0\pmod p}}\Bigl
(\frac
12-\frac{\nu}p\Bigr )\log p\label{T1}\\
T_2(n)&=\sum_{\substack{0\le \nu <p<Cn\\f(\nu)\equiv 0\pmod p}}\Bigl
(\frac 12-\Bigl \{\frac{n-\nu}p\Bigr \}\Bigr )\log p.\label{T2}
\end{align}
Sums $T_1(n)$ and $T_2(n)$ will be $o(n)$ as a consequence of
theorem \ref{I}. But it is not completely obvious and we will do it
in detail in the next subsection.

Now we will simplify \eqref{Ln4} a little more in the next lemma.
\begin{lem}
\begin{equation}\label{Ln5}
\log L_n(f)=n\log n+cn+S_2(n)-T_1(n)-T_2(n)+O(n^{2/3}\log n),
\end{equation}
where
\begin{align*}
c=\log a-\log C-2+\gamma-\sum_{p\nmid 2aD}\frac{(d/p)\log
p}{p-1}+\sum_{p\mid 2aD}\log p\Bigl(\frac 1{p-1}-\sum_{k\ge
1}\frac{s(f,p^k)}{p^k}\Bigr )\notag
\end{align*}
and $S_2(n),\ T_1(n)$ and $T_2(n)$ are as in \eqref{S2}, \eqref{T1}
and \eqref{T2}.
\end{lem}
\begin{proof} Let $D=l^2d$ where $d$ is a fundamental discriminant. First we observe that $(D/p)=(l/p)^2(d/p)$ and that if $p\nmid D$ then $(D/p)=(d/p)$.
As a consequence of the prime number theorem on arithmetic
progressions we know that the sum  $\sum_{p}\frac{(d/p)\log p}{p-1}$
is convergent. On the other hand, the well known estimate
$\sum_{p\le x}\frac{\log p}{p-1}=\log x-\gamma+o(1)$ where $\gamma $
is the Euler constant, implies  that
\begin{align}\label{SSS}
\sum_{\substack{p< Cn\\ p\nmid 2aD}}\frac{(1+(D/p))\log
p}{p-1}&=\log n+\log C -\gamma-\sum_{p\mid 2aD}\frac{\log p}{p-1}
\\&+\sum_{p\nmid 2aD}\frac{(d/p)\log p}{p-1}+o(1).\notag
\end{align}
 Finally we put \eqref{SSS}  in \eqref{Ln4}.
\end{proof}

\subsection{Equidistribution of the roots $\pmod p$ of a quadratic polynomial}
Now we develop a technology to prove that $T_1(n)$, $T_2(n)$ and
other similar sums which will appear in the estimate of $S_2(n)$ are
all $o(n)$.

These sums are all of the form
\begin{equation}\label{form}
\sum_{\substack{0\le \nu<p \le x,\ p\in S\\
f(\nu)\equiv 0\pmod p}}a(\nu,p,x)\log p
\end{equation}
for some function $a(\nu,p,x)\ll 1$. By partial summation we also
get easily that
\begin{align}\label{cerolog} \sum_{\substack{0\le \nu<p \le x,\ p\in S\\
f(\nu)\equiv 0\pmod p}}a(\nu,p,x)\log p&=\log x\sum_{\substack{0\le \nu<p \le x,\ p\in S\\
f(\nu)\equiv 0\pmod p}}a(\nu,p,x)-\int_1^x\frac 1t \sum_{\substack{0\le \nu<p \le t,\ p\in S\\
f(\nu)\equiv 0\pmod p}}a(\nu,p,x)\\&=\log x\sum_{\substack{0\le \nu<p \le x,\ p\in S\\
f(\nu)\equiv 0\pmod p}}a(\nu,p,x)+o(x/\log x).\end{align} Hence, to
prove that the sums \eqref{form} are $o(x)$ we must prove that $$\sum_{\substack{0\le \nu<p \le x,\ p\in S\\
f(\nu)\equiv 0\pmod p}}a(\nu,p,x)=o(x/\log x).$$

Theorem \ref{I} implies, in particular, that for any arithmetic
progression $S$  and for any piecewise continuos function $g$ in
$[0,1]$ such that $\int_0^1g=0$ we
have that \begin{equation}\label{cero} \sum_{\substack{0\le \nu<p \le x,\ p\in S\\
f(\nu)\equiv 0\pmod p}}g(\nu/p)=o(x/\log x).\end{equation}

\begin{lem}\label{S1} Let $f$ be an irreducible polynomial in $\Z[x]$.
We have that the sums $T_1(n)$ and $T_2(n)$ defined in \eqref{T1}
and \eqref{T2} are both $o(n)$.
\end{lem}
\begin{proof}
To prove that $T_1(n)=o(n)$  we apply \eqref{cero} to the function
$g(x)=x-1/2$.

To prove that $T_2(n)=o(n)$ the strategy  is to split the range of
the primes in small intervals such that  $n/p$ are almost constant
in each interval. We take $H$ a large, but a fixed number and we
divide the interval $[1,Cn]$ in $H$ intervals
$L_h=(\frac{h-1}HCn,\frac hHCn],\ h=1,\dots , H$. Now we write
\begin{equation}\sum_{\substack{0\le \nu <p<n\\ f(\nu)\equiv 0\pmod p}}\Bigl(\Bigl \{\frac{n-\nu}p\Bigr\}-\frac 12\Bigr )=
\Sigma_{31}+\Sigma_{32}+\Sigma_{33}+O(n/(H^{1/3}\log n))
\end{equation}
 where
\begin{align*}
\Sigma_{31}&=\sum_{H^{2/3}\le h\le H}\sum_{\substack{0\le \nu<p\in L_h\\
f(\nu)\equiv 0\pmod p}}\Bigl(  \Bigl\{ \frac{H}h-\frac{\nu}p \Bigr \}-\frac 12  \Bigr )\\
\Sigma_{32}&=\sum_{H^{2/3}\le h\le H}\sum_{\substack{0\le \nu<p\in L_h\\
f(\nu)\equiv 0\pmod p\\ \frac{\nu}p\not \in [\frac Hh,\frac
H{h-1}]}}\Bigl (\Bigl \{ \frac{n}p-\frac{\nu}p\Bigr \}-  \Bigr \{
\frac{H}h-\frac{\nu}p \Bigr \}  \Bigr )\\
\Sigma_{33}&=\sum_{H^{2/3}\le h\le H}\sum_{\substack{0\le \nu<p\in L_h\\
f(\nu)\equiv 0\pmod p\\ \frac{\nu}p \in [\frac Hh,\frac
H{h-1}]}}\Bigl (\Bigl \{ \frac{n}p-\frac{\nu}p\Bigr \}-  \Bigl \{
\frac{H}h-\frac{\nu}p \Bigr \}  \Bigr ).
\end{align*}
To estimate $\Sigma_{31}$ we apply \eqref{cero} with the function
$\Bigl\{ \frac{H}h-x \Bigr \}-\frac 12$ in each $L_h$ and we obtain
\begin{equation}\label{31}\Sigma_{31}=o(Hn/\log n)=o(n\log n)
\end{equation}since $H$ is a constant.

To bound $\Sigma_{32}$ we observe that if $p\in L_h$ and
$\frac{\nu}p\not \in [\frac Hh,\frac H{h-1}]$, then
$$0\le \Bigl \{ \frac{n}p-\frac{\nu}p\Bigr \}- \Bigl \{
\frac{H}h-\frac{\nu}p \Bigr \}=\frac{n}p-\frac{H}h\le \frac
H{h(h-1)}.$$ Thus
\begin{align}\label{32}
|\Sigma_{32}|\ll \sum_{H^{2/3}\le h<H}\sum_{p\in L_h}\frac H{h^2}\ll
\sum_{H^{2/3}\le h<H}\sum_{p\in L_h}\frac{1}{H^{1/3}}\ll
\frac{\pi(n)}{H^{1/3}} \ll \frac{n}{H^{1/3}\log n}.
\end{align}

To bound $\Sigma_{33}$ first we observe that
\begin{align*}\Sigma_{33}&\ll  \sum_{H^{2/3}\le h<H}\sum_{\substack{0\le \nu<p\in L_h\\
f(\nu)\equiv 0\pmod p\\ \frac{\nu}p \in [\frac Hh,\frac
H{h-1}]}}1\\ &=\sum_{H^{2/3}\le h<H}\sum_{\substack{0\le \nu<p\in L_h\\
f(\nu)\equiv 0\pmod p}}\Bigl (\chi_{[H/h,H/(h-1)]}(\nu/p)-\frac
H{h(h-1)}\Bigr )\\ &+\sum_{H^{2/3}\le h<H}\sum_{\substack{0\le \nu<p\in L_h\\
f(\nu)\equiv 0\pmod p}}\frac H{h(h-1)},\end{align*} where, here and
later, $\chi_{[a,b]}(x)$ denotes the characteristic function of the
interval $[a,b]$.

Theorem \ref{I} implies that $$\sum_{\substack{0\le \nu<p\in L_h\\
f(\nu)\equiv 0\pmod p}}\Bigl (\chi_{[H/h,H/(h-1)]}(\nu/p)-\frac
H{h(h-1)}\Bigr )=o(n/\log n).$$ Thus,
\begin{eqnarray}\label{33}\Sigma_{33}&\ll &\sum_{H^{2/3}\le h<H} o(\frac n{\log n})+\sum_{H^{2/3}\le h<H}\sum_{\substack{0\le \nu<p\in L_h\\
f(\nu)\equiv 0\pmod p}}\frac 1{H^{1/3}}\\ &\ll& o(n/\log
n)+\frac{\pi(n)}{H^{1/3}}\ll o(n/\log n)+O(n/H^{1/3}\log
n).\notag\end{eqnarray}

Estimates \eqref{31}, \eqref{32} and \eqref{33} imply $\Sigma_3\ll
o(n/\log n)+n/(H^{1/3}\log n)$. Since $H$ can be chosen arbitrarily
large we have that $\Sigma_3=o(n/\log n)$ which finish the proof.
\end{proof}

To present lemma \ref{K} we need some preparation.

For primes $p\in \mathcal P_f$  the congruence $f(x)\equiv 0\pmod p$
has exactly two solutions, say $0\le \nu_{p,1},\ \nu_{p,2}<p$.

In some parts of the proof of theorem \ref{principal} we will need
estimate some quantities depending on $\min (\nu_{p,1},\nu_{p,2})$.
For this reason it is convenient to know how they are related.

If $f(x)=ax^2+bx+c$ and $p\in \mathcal P_f$ then
$\nu_{p,1}+\nu_{p,2}\equiv -b/a \pmod p$. Next lemma will give more
information when the prime $p$ belongs to some particular arithmetic
progression.
\begin{lem}\label{1+2}Let $q=a/(a,b),\
l=b/(a,b)$. For any $r,\  (r,q)=1$ and for any prime $p\equiv
lr^{-1}\pmod q$ and $p\in \mathcal P_f$ we have
\begin{equation}
\frac{\nu_{p,1}}p+\frac{\nu_{p,2}}p\equiv \frac rq-\frac{l}{pq}\pmod
1
\end{equation}
\end{lem}
\begin{proof}
To avoid confusions we denote by $\overline q_p$ and $\overline p_q$
the inverses of $q\pmod p$ and $p\pmod q$ respectively. From the
obvious congruence $q\overline q_p+p\overline p_q\equiv 1\pmod{pq}$
we deduce that $\frac{\overline q_p}p+\frac{\overline p_q}q-\frac
1{pq}\in \Z$. Since $p\equiv l\overline r_q\pmod q$ we obtain
$\frac{\overline q_p}p\equiv \frac 1{pq}-\frac{r\overline l_q}q\pmod
1$. Thus
$$\frac{\nu_{p,1}}p+\frac{\nu_{p,2}}p\equiv \frac{-l\overline
q_p}p\equiv -l\left (\frac 1{pq}-\frac{r\overline l_q}q\right
)\equiv \frac rq-\frac l{pq}\pmod 1.$$
\end{proof}
Since the two roots are symmetric respect to $\frac
r{2q}-\frac{l}{2pq}$, necessarily one of then lies in $\Bigl [\frac
r{2q}-\frac{l}{2pq}, \frac 12+\frac r{2q}-\frac{l}{2pq}\Bigr )\pmod
1$ and the other in the complementary set.
\begin{defn}For $(r,q)=1,\ 1\le r\le q$, $p\equiv lr^{-1}\pmod q$ and $p\in \mathcal P_f$ we define $\nu_{p,1}$ the root
of $f(x)\equiv \pmod p$ such that
$$\frac{\nu_{p,1}}p\in T_{rp}=\Bigl[\frac r{2q}-\frac{l}{2pq},
\frac 12+\frac r{2q}-\frac{l}{2pq}\Bigr )\pmod 1,$$ and we define
$\nu_{p,2}$ the root of $f(x)\equiv 0\pmod p$ such that
$\frac{\nu_{p,2}}p\in [0,1)\setminus T_{rp}$.
\end{defn}

\begin{lem}\label{K}Assume the notation above. Let
$\alpha_1,\alpha_2,\beta_1,\beta_2,c_1,c_2$ be constants  and
$g_1(x),g_2(x)$ two linear functions satisfying that
$$J_n(p)=\left [g_1\left (\frac{n}p\right )+\frac{c_1}p,g_2\left (\frac np\right )+
\frac{c_2}p\right ]\subset T_{rp}$$ for any prime $p\in
K_n=[\alpha_1n+\beta_1,\alpha_2n+\beta_2]$. We have
\begin{equation}\label{sum}
\sum_{\substack{p\in K_n\cap  \mathcal P_f\\ p\equiv lr^{-1}\pmod
q}}\Bigl (\chi_{J_n(p)}\Bigl (\frac{\nu_{p,1}}p\Bigr
)-2|J_n(p)|\Bigr )\log p=o(n)
\end{equation}
where $\chi_I$ is the characteristic function of the set $I$.
\end{lem}
\begin{proof}
Since $J_n(p)\subset T_{rp}$ then $\nu_2/p\not \in J_n(p)$ and we
can write
$$\sum_{\substack{p\in K_n\cap \mathcal P_f\\ p\equiv lr^{-1}\pmod
q}}\chi_{J_n(p)}\Bigl(\frac{\nu_{p,1}}p\Bigr )\log
p=\sum_{\substack{1\le \nu\le p\in K_n,\\ f(\nu)\equiv 0\pmod
p\\p\equiv lr^{-1}\pmod q}}\chi_{J_n(p)}\Bigl (\frac{\nu}p\Bigr
)\log p
$$
and $$\sum_{\substack{p\in K_n\cap \mathcal P_f\\ p\equiv
lr^{-1}\pmod q}}2|J_n(p)|\log p=\sum_{\substack{1\le \nu\le p\in K_n,\\
f(\nu)\equiv 0\pmod p\\p\equiv lr^{-1}\pmod q}}|J_n(p)|\log p.$$
Thus,
$$\sum_{\substack{p\in K_n\cap \mathcal P_f\\ p\equiv lr^{-1}\pmod
q}}\Bigl (\chi_{J_n(p)}\Bigl(\frac{\nu_{p,1}}p\Bigr)-2|J_n(p)|\Bigr
)\log p=\sum_{\substack{1\le \nu\le p\in K_n,\\ f(\nu)\equiv 0\pmod
p\\p\equiv lr^{-1}\pmod q}}\Bigl (\chi_{J_n(p)}\Bigl
(\frac{\nu}p\Bigr )-|J_n(p)|\Bigr )\log p. $$ To prove the lemma is
enough to prove that
\begin{equation}
\sum_{\substack{1\le \nu\le p\in K_n,\\ f(\nu)\equiv 0\pmod
p\\p\equiv lr^{-1}\pmod q}}\Bigl (\chi_{J_n(p)}\Bigl
(\frac{\nu}p\Bigr )-|J_n(p)|\Bigr )=o(n/\log n).
\end{equation}

We proceed as above. We split $K_n$ in intervals
$L_h=(\frac{h-1}Hn,\frac hHn]$ of length $n/H$ and two extra
intervals $I,\ F$ (the initial and the final intervals) of length
$\le n/H$. Here $h$ runs over a suitable set of consecutive integers
$\mathcal H$ of cardinality $\ll (\alpha_2-\alpha_1)H$.

Let $I_h$ denote the interval
$[g_1(H/h)+c_1H/(nh),g_2(H/h)+c_2H/(nh)]$.

We write
\begin{equation}
\sum_{\substack{1\le \nu\le p\in K_n,\\ f(\nu)\equiv 0\pmod
p\\p\equiv lr^{-1}\pmod q}}\Bigl
(\chi_{J_n(p)}\Bigl(\frac{\nu}p\Bigr )-|J_n(p)|\Bigr
)=\Sigma_1+\Sigma_2+\Sigma_3+\Sigma_4
\end{equation}
where
\begin{align*}
\Sigma_1&=\sum_{h\in \mathcal H}\sum_{\substack{0\le \nu<p\in
L_h\\f(\nu)\equiv 0\pmod p\\ p\equiv lr^{-1}\pmod q}}\Bigl
(\chi_{I_h}\Bigl (\frac{\nu}p\Bigr)-|I_h|\Bigr)\\
\Sigma_2&=\sum_{h\in \mathcal H}\sum_{\substack{0\le \nu<p\in
L_h\\f(\nu)\equiv 0\pmod p\\
p\equiv lr^{-1}\pmod q}}(|I_h|-|J_n(p)|)\\
\Sigma_3&=\sum_{h\in \mathcal H}\sum_{\substack{0\le \nu<p \in
L_h\\f(\nu)\equiv 0\pmod p\\ p\equiv lr^{-1}\pmod q}}\Bigl
(\chi_{J_n(p)}\Bigl(\frac{\nu}p\Bigr )-\chi_{I_h}\Bigl
(\frac{\nu}p\Bigr )\Bigr )\\
\Sigma_4&=\sum_{\substack{0\le \nu<p\in I\cup F\\f(\nu)\equiv 0\pmod
p\\ p\equiv lr^{-1}\pmod q}}\Bigl(\chi_{I_h}\Bigl (\frac{\nu}p\Bigl
)-|J_n(p)|\Bigr ).
\end{align*}

The inner sum in $\Sigma_1$ can be estimated as we did in lemma
\ref{S1}, (with the function $g(x)=\chi_I(x)-|I|$ instead of
$g(x)=x-1/2$), and we get again that $\Sigma_1=o(n/\log n)$.

To estimate $\Sigma_2$ and $\Sigma_3$ we observe that if $p\in L_h$
then $J_n(p)$ and $I_h$ are almost equal. Actually, comparing the
end points of both intervals and because $g$ is a linear function,
we have $|J_n(p)|-|I_h|\ll \min(1,H/h^2)$ and
$\chi_{J_n(p)}(x)=\chi_{I_h}(x)$ except for a set $E_h$ of measure
$\ll \min (1,H/h^2)$.

We have
\begin{align*}
\Sigma_2&\ll  \sum_{h\in \mathcal H} \sum_{p\in L_n}\min(1,H/h^2)
\ll \sum_{h\le H^{2/3}}\sum_{p\in L_h}+\sum_{H^{2/3}<h\in \mathcal H}\sum_{p\in L_h}\frac 1{H^{1/3}}\\
&\ll\pi(n/H^{1/3})+\frac 1{H^{1/3}}\pi(\alpha_1n+\alpha_2)\ll
n/(H^{1/3}\log n).
\end{align*}
To bound $\Sigma_{3}$ first we observe that
 \begin{align*} \Sigma_3&\ll \sum_{h\in \mathcal
H}\sum_{\substack{0\le \nu<p\in L_h\\f(\nu)\equiv 0\pmod p\\ p\equiv
lr^{-1}\pmod
q}}\chi_{E_h}(\nu/p)\\
&=\sum_{h\in \mathcal H}\sum_{\substack{0\le \nu<p\in
L_h\\f(\nu)\equiv 0\pmod p\\ p\equiv lr^{-1}\pmod q}}\left
(\chi_{E_h}(\nu/p)-|E_h|\right )+\sum_{h\in \mathcal
H}\sum_{\substack{0\le \nu<p\in L_h\\f(\nu)\equiv 0\pmod p\\ p\equiv
lr^{-1}\pmod q}}|E_h|.\end{align*} Theorem \ref{I} implies that
$$\sum_{\substack{0\le \nu<p\in
L_h\\f(\nu)\equiv 0\pmod p\\ p\equiv lr^{-1}\pmod q}}\left
(\chi_{E_h}(\nu/p)-|E_h|\right )=o(n/\log n).  $$ On the other hand,
\begin{align*}\sum_{h\in \mathcal
H}\sum_{\substack{0\le \nu<p\in L_h\\f(\nu)\equiv 0\pmod p\\ p\equiv
lr^{-1}\pmod q}}|E_h|&\ll \sum_{h\le H^{2/3}}\sum_{p\in
L_h}1+\sum_{H^{2/3}<h\in \mathcal H}\sum_{p\in L_h}\frac H{h^2}\\
&\ll \pi(n/H^{1/3})+\frac 1{H^{1/3}}\pi(\alpha_1n+\alpha_2)\\
&\ll \frac n{H^{1/3}\log n}.
\end{align*}
Thus, $\Sigma_3\ll o(n/\log n)+n/(H^{1/3}\log n)$.

Finally we estimate $\Sigma_4$. We observe that $$|\Sigma_4|\le
\sum_{p\in I}1+ \sum_{p\in F}1\ll n/(H\log n) $$ as a consequence of
the prime number theorem. Then
$$\Sigma_1+\Sigma_2+\Sigma_3+\Sigma_4 =O(n/(H^{1/3}\log n))+O(n/(H\log n))+o(n/\log n)$$
which finish the proof because we can take $H$ arbitrarily large.
\end{proof}

\subsection{Estimate of $S_2(n)$ and end of the proof}
\begin{lem}
\begin{equation}\label{S22}S_2(n)=n\Bigl (1+\log C-\log 4+\frac
1{\phi(q)}\sum_{(r,q)=1}\log(1+\frac rq)\Bigr )+o(n)\end{equation}
\end{lem}
\begin{proof}
Following the notation of lemma \ref{1+2} we split
$$S_2(n)=\sum_{\substack{(r,q)=1\\ 1\le r\le q}}S_{2r}(n)+\sum_{p\le
l}\beta_p^*(n)\log p=\sum_{\substack{(r,q)=1\\ 1\le r\le
q}}S_{2r}(n) +O(1)$$ where
\begin{equation}
S_{2r}(n)=\sum_{\substack{l<p\le Cn\\p\equiv lr^{-1}\pmod
q}}\beta_p^*(n)\log p.
\end{equation}

 Since $p\equiv lr^{-1}\pmod q$, lemma
\ref{1+2} implies that $\frac{\nu_{p,1}}p+\frac{\nu_{p,2}}p\equiv
\frac rq-\frac l{pq}\pmod 1$. We observe also that, since $p>l$ we
have that $0< \frac rq-\frac l{pq}\le 1$.

Now we will check that
$$\beta_p^*(n)=\begin{cases}1,\ &\text{ if } \qquad \qquad \quad \frac np\ge \frac 12+\frac
r{2q}-\frac{l}{2pq}\\
\chi_{[\frac r{2q}-\frac{l}{2pq},\frac np]}(\nu_{p,1}/p),\ &\text{
if }\quad  \frac rq-\frac{l}{pq}<\frac np< \frac 12+\frac
r{2q}-\frac{l}{2pq}\\
\chi_{[\frac r{2q}-\frac{l}{2pq},\frac
r{q}-\frac{l}{pq}]}(\nu_{p,1}/p),\ &\text{ if }\ \frac
r{2q}-\frac{l}{2pq}\le \frac np \le \frac rq-\frac{l}{pq}\\
\chi_{[\frac rq-\frac{l}{pq}-\frac np,\frac
rq-\frac{l}{pq}]}(\nu_{p,1}/p)&\text{ if }\qquad \qquad \quad \
\frac np<\frac r{2q}-\frac{l}{2pq}
\end{cases}$$
We observe that $\beta_p^*(n)=1$ if and only if
$\frac{\nu_{p,1}}p\le \frac np$ or $\frac{\nu_{p,2}}p\le \frac np$.
We remind that
\begin{equation}\label{menor}\frac r{2q}-\frac l{2pq}\le
\frac{\nu_{p,1}}p< \frac 12+\frac r{2q}-\frac l{2pq}\end{equation}
Also we observe that lemma \ref{1+2} implies that
\begin{equation}\label{igual}\frac{\nu_{p,2}}p=\begin{cases}\frac r{q}-\frac l{pq}-\frac{\nu_{p,1}}p &\text{  if }\frac{\nu_{p,1}}p\le \frac r{q}-\frac l{pq}
\\\frac r{q}-\frac l{pq}-\frac{\nu_{p,1}}p+1&\text{  if }\frac{\nu_{p,1}}p> \frac r{q}-\frac l{pq}.\end{cases}\end{equation}

\begin{itemize}
    \item Assume $\frac np\ge \frac 12+\frac
r{2q}-\frac{l}{2pq}$. Then $ \nu_{p,1}< p\left (\frac 12+\frac
r{2q}-\frac l{2pq}\right )<n$, so $\beta_p^*(n)=1$
    \item Assume $\frac rq-\frac{l}{pq}<\frac np< \frac 12+\frac
r{2q}-\frac{l}{2pq}$.
\begin{itemize}
    \item If $\chi_{[\frac r{2q}-\frac l{2pq}, \frac np]}(\nu_{p,1}/p)=1$
    then
    $\nu_{p,1}\le n,$ so $\beta_p^*(n)=1$.
    \item If $\chi_{[\frac r{2q}-\frac l{2pq}, \frac np]}(\nu_{p,1}/p)=0$
    then  $\frac{\nu_{p,1}}p>\frac np>\frac rq-\frac l{pq}$.
    Relations (\ref{menor}) and (\ref{igual})
    imply that
    $\frac{\nu_{p,2}}p=1+\frac rq-\frac l{pq}-\frac{\nu_{p,1}}p>\frac 12+\frac r{2q}-\frac l{2pq}>\frac
    np.$
    Since $\nu_{p,1}>n$ and $\nu_{p,2}>n$ we get $\beta_p^*(n)=0$.
\end{itemize}
    \item Assume $\frac
r{2q}-\frac{l}{2pq}\le \frac np \le \frac rq-\frac{l}{pq}.$
\begin{itemize}
    \item If $\chi_{[\frac{r}{2q}-\frac l{2pq},\frac{r}{q}-\frac
    l{pq}]}(\nu_{p,1}/p)=1$ then (\ref{igual}) imply that $0<
    \frac{\nu_{p,2}}p\le \frac{r}{2q}-\frac l{2pq}$, which implies
    that $\nu_{p,2}\le n$, so $\beta_p^*(n)=1$.
    \item If $\chi_{[\frac{r}{2q}-\frac l{2pq},\frac{r}{q}-\frac
    l{pq}]}(\nu_{p,1}/p)=0$ then $\frac{\nu_{p,1}}p>\frac{r}{q}-\frac
    l{pq}\ge \frac np$ and relation (\ref{igual}) imply that
    $\frac{\nu_{p,2}}p=\frac r{q}-\frac l{pq}-\frac{\nu_{p,1}}p+1>\frac r{q}-\frac l{pq}\ge \frac
    np.$ Since $\nu_{p,1}>n$ and $\nu_{p,2}>n$ we get $\beta_p^*(n)=0$.
\end{itemize}
    \item Assume $\frac np<\frac r{2q}-\frac l{2pq}$.
    \begin{itemize}
        \item If $\chi_{[\frac rq-\frac{l}{pq}-\frac np,\frac
rq-\frac{l}{pq}]}(\nu_{p,1}/p)=1$ then $\frac{\nu_{p,1}}p\le \frac
rq-\frac l{pq}$ and relation (\ref{igual}) implies that
$\frac{\nu_{p,2}}p=\frac rq-\frac l{pq}-\frac{\nu_{p,1}}p\le \frac
rq-\frac l{pq}-(\frac rq-\frac{l}{pq}-\frac np )=\frac np,$ so
$\beta_p^*(n)=1$
        \item If $\chi_{[\frac rq-\frac{l}{pq}-\frac np,\frac
rq-\frac{l}{pq}]}(\nu_{p,1}/p)=0$ we distinguish two
cases:\begin{itemize}
    \item If $\frac r{2q}-\frac{l}{2q}\le \frac{\nu_{p,1}}p<\frac
    rq-\frac l{pq}-\frac np$ then $\frac{\nu_1,p}p\ge \frac
    r{2q}-\frac l{2q}>\frac np$, and also we have that
    $\frac{\nu_{p,2}}p=\frac rq-\frac l{pq}-\frac{\nu_{p,1}}p>\frac rq-\frac
    l{pq}-\left( \frac
    rq-\frac l{pq}-\frac np \right )=\frac np$. Thus $\beta_p^*(n)=0$
    \item If $\frac r{q}-\frac{l}{pq}<\frac{\nu_{p,1}}p<
   \frac 12 +\frac r{2q}-\frac l{2pq}$ then $\frac{\nu_{p,1}}p>\frac
   12\left (\frac r{q}-\frac{l}{pq}\right )>\frac np$. On the other hand,
   $\frac{\nu_{p,2}}p=\frac rq-\frac l{pq}-\frac{\nu_{p,1}}p+1>\frac rq-\frac l{pq}-\left ( \frac 12 +\frac r{2q}-\frac l{2pq}  \right )+1=
   \frac 12+\frac r{2q}-\frac l{2pq}>\frac np.$ Thus, again we have
   that $\beta_p^*(n)=0$.
\end{itemize}
    \end{itemize}
\end{itemize}
Now we split $S_{2r}(n)=\sum_{i=1}^4S_{2ri}(n)$  according the
ranges of the primes involved in lemma above.
\begin{eqnarray*}
S_{2r1}(n)&=&\sum_{\substack{l<p\le
\frac{n+l/(2q)}{1/2+r/(2q)}\\p\equiv
lr^{-1}\pmod q\\ p\in \mathcal P_f}}\log p\\
S_{2r2}(n)&=&\sum_{\substack{\frac{n+l/(2q)}{1/2+r/(2q)}<p<
\frac{n+l/q}{r/q}\\p\equiv lr^{-1}\pmod q\\ p\in \mathcal
P_f}}\chi_{[\frac r{2q}-\frac
l{2pq},\frac np]}(\nu_{p,1}/p)\log p\\
S_{2r3}(n)&=&\sum_{\substack{\frac qr(n+\frac{l}q)\le p\le
\frac{2q}r(n+\frac{l}q)\\p\equiv lr^{-1}\pmod q\\p\in \mathcal
P_f}}\chi_{[\frac r{2q}-\frac{l}{2pq},\frac
r{q}-\frac{l}{pq}]}(\nu_{p,1}/p)\log p\\
S_{2r4}(n)&=&\sum_{\substack{\frac{2q}r(n+\frac l{2q})<p<Cn\\p\equiv
lr^{-1}\pmod q\\p\in\mathcal P_f}}\chi_{[\frac rq-\frac{l}{pq}-\frac
np,\frac rq-\frac{l}{pq}]}(\nu_{p,1}/p)\log p.
\end{eqnarray*}

Since $(q,D)=1$ and the primes are odd numbers, the primes $p\equiv
lr^{-1}\pmod q,\ p\in \mathcal P_f$ lie in a set of
$\phi(4qD)/(2\phi(q))$ arithmetic progressions modulo $4qD$. The
prime number theorem for arithmetic progressions implies that
\begin{equation}\label{pa}\sum_{\substack{p\le x\\p\equiv lr^{-1}\pmod q,\ p\in
P_f}}\log p\sim \frac{x}{2\phi(q)}\end{equation} and
\begin{equation}\label{pa2}\sum_{\substack{ax<p\le bx\\p\equiv lr^{-1}\pmod q,\ p\in P_f}}\frac{\log p}p= \frac{\log
(b/a)}{2\phi(q)}+o(1)\end{equation}

We will use these estimates and lema \ref{K} to estimate
$S_{2ri}(n),\ i=1,2,3,4.$

By (\ref{pa}) we have
\begin{equation}
S_{2r1}(n)=\frac n{\phi(q)}\frac q{q+r}+o(n).
\end{equation}
To estimate $S_{5r2}$ we write
\begin{eqnarray*}
S_{2r2}(n)&=&\sum_{\substack{\frac{n+l/(2q)}{1/2+r/(2q)}<p<
\frac{n+l/q}{r/q}\\p\equiv lr^{-1}\pmod q}}\chi_{[\frac r{2q}-\frac
l{2pq},\frac np]}(\nu_{p,1}/p)\log p\\
&=&\sum_{\substack{\frac{n+l/(2q)}{1/2+r/(2q)}<p<
\frac{n+l/q}{r/q}\\p\equiv lr^{-1}\pmod q}}\Bigl (\frac{2n}p- \frac
r{q}+\frac l{pq}\Bigr )\log p\\
&+&\sum_{\substack{\frac{n+l/(2q)}{1/2+r/(2q)}<p<
\frac{n+l/q}{r/q}\\p\equiv lr^{-1}\pmod q}}\Bigl (\chi_{[\frac
r{2q}-\frac l{2pq},\frac np]}(\nu_{p,1}/p)-2\Bigl (\frac np-\frac
r{2q}+\frac l{2pq}\Bigr )\Bigr )\log p
\end{eqnarray*}

Lemma \ref{K} implies that the last sum is $o(n)$.  Thus,
\begin{eqnarray*}
S_{5r2}&=&\sum_{\substack{\frac{n+l/(2q)}{1/2+r/(2q)}<p<
\frac{n+l/q}{r/q}\\p\equiv lr^{-1}\pmod q\\ p\in \mathcal P_f}}\Bigl
(\frac{2n}p- \frac
r{q}\Bigr )\log p   +o(n)\\
&=&2n\sum_{\substack{\frac{n+l/(2q)}{1/2+r/(2q)}<p<
\frac{n+l/q}{r/q}\\p\equiv lr^{-1}\pmod q\\ p\in \mathcal
P_f}}\frac{\log p}p-\frac rq
\sum_{\substack{\frac{n+l/(2q)}{1/2+r/(2q)}<p<
\frac{n+l/q}{r/q}\\p\equiv lr^{-1}\pmod q\\ p\in \mathcal P_f}}\log p +o(n)\\
&=&\frac n{\phi(q)}\log\Bigl(\frac 12+\frac q{2r}\Bigr )-\frac
n{\phi(q)}\Bigl (\frac 12-\frac r{q+r}\Bigr )+o(n)
\end{eqnarray*}
by (\ref{pa}) and (\ref{pa2}).

To estimate $S_{2r3}(n)$ we write
\begin{eqnarray*}S_{2r3}(n)&=&\sum_{\substack{\frac qr(n+\frac{l}q)\le p\le
\frac{2q}r(n+\frac{l}q)\\p\equiv lr^{-1}\pmod q\\p\in \mathcal
P_f}}\Bigl (\frac r{q}-\frac{l}{pq}\Bigr)\log p\\
&+&\sum_{\substack{\frac qr(n+\frac{l}q)\le p\le
\frac{2q}r(n+\frac{l}q)\\p\equiv lr^{-1}\pmod q\\p\in \mathcal
P_f}}\Bigl (\chi_{[\frac r{2q}-\frac{l}{2pq},\frac
r{q}-\frac{l}{pq}]}(\nu_{p,1}/p)- \Bigl (\frac
r{q}-\frac{l}{pq}\Bigr )\Bigr )\log p\\
&=&\frac{n}{2\phi(q)}+o(n)
\end{eqnarray*}
by (\ref{pa}) and lema \ref{K}.

To estimate $S_{2r4}(n)$ we write
\begin{eqnarray*}
S_{2r4}(n)&=&\sum_{\substack{\frac{2q}r(n+\frac l{2q})<p<Cn\\p\equiv
lr^{-1}\pmod q\\p\in\mathcal P_f}}\Bigl (\frac{2n}p+\frac{2l}{pq}\Bigr )\log p\\
&=&\sum_{\substack{\frac{2q}r(n+\frac l{2q})<p<Cn\\p\equiv
lr^{-1}\pmod q\\p\in\mathcal P_f}}\Bigl (\chi_{[\frac
rq-\frac{l}{pq}-\frac np,\frac rq-\frac{l}{pq}]}(\nu_{p,1}/p)-\Bigl
(\frac{2n}p+\frac{2l}{pq}\Bigr ) \Bigr)\log p\\
&=&\frac n{\phi(q)}\Bigl(\log C-\log(2q/r)\Bigr )+o(n)
\end{eqnarray*}
by (\ref{pa2}) and lemma \ref{K}.

Thus
\begin{eqnarray*}S_{2r}(n)&=&S_{2r1}(n)+S_{2r2}(n)+S_{2r3}(n)+S_{2r4}(n)+O(1)\\&=&\frac
n{\phi(q)}\frac q{q+r}+o(n)\\ &+&\frac n{\phi(q)}\log \Bigl(\frac
12+\frac q{2r}\Bigr )-\frac n{\phi(q)}\Bigl (\frac 12-\frac
r{q+r}\Bigr )+o(n)\\ &+&\frac{n}{2\phi(q)}+o(n)\\&+&\frac
n{\phi(q)}\Bigl (\log C-\log(2q/r)\Bigr )+o(n)\\&=&\frac
n{\phi(q)}\Bigl (1+ \log C-\log 4+\log(1+r/q)\Bigr
)+o(n).\end{eqnarray*}

 Now sum in all $r\le q,\ (r,q)=1$ to finish the estimate of
$S_2(n)$.
\end{proof}

Finally we substitute \eqref{S22} in \eqref{Ln5} to conclude the
proof of theorem \ref{principal}.

\section{Computation of the constant $B_f$} %
The sum $\sum_{p}\frac{(d/p)\log p}{p-1}$, appearing in the formula
of the constant $B_f$ converges very slowly. Next lemma gives an
alternative expression for this sum, more convenient to obtain a
fast computation.
\begin{lem}
\begin{align}
\sum_{p}\frac{(d/p)\log
p}{p-1}=\sum_{k=1}^{\infty}\frac{\zeta'(2^k)}{\zeta(2^k)}
-\sum_{k=0}^{\infty}\frac{L'(2^k,\chi_d)}{L(2^k,\chi_d)}+\sum_{p\mid
d}s_p.
\end{align}
where $s_p=\sum_{k=1}^{\infty}\frac{\log p}{p^{2^k}-1}$.
\end{lem}
\begin{proof}For $s>1$ we consider the function $G_d(s)=\prod_p\Bigl (1-\frac
1{p^s}\Bigr )^{(d/p)}$. Taking the derivative of the logarithm of
$G_d(s)$ we obtain that
\begin{equation}\label{G}\frac{G_d'(s)}{G_d(s)}=\sum_p\frac{(d/p)\log p}{p^s-1}.\end{equation}

Since $L(s,\chi_d)=\prod_p\Bigl (1-\frac{(d/p)}p^s\Bigr )^{-1}$ we
have
\begin{align}G_d(s)L(s,\chi_d)&=\prod_p\Bigl (1-\frac 1{p^s}\Bigr
)^{(d/p)}\Bigl (1-\frac{(d/p)}{p^s}\Bigr
)^{-1}\\&=\prod_{(d/p)=-1}\Bigl (1-\frac 1{p^{2s}}\Bigr )^{-1}\\
&=\prod_p \Bigl (1-\frac 1{p^{2s}}\Bigr
)^{\frac{(d/p)-1}2}\prod_{p\mid d}\Bigl (1-\frac 1{p^{2s}}\Bigr
)^{1/2}\\&=G_d^{1/2}(2s)\zeta^{1/2}(2s)T^{1/2}(2s)\end{align} where
$T(s)=\prod_{p\mid d}\Bigl (1-\frac 1{p^s}\Bigr)$.

The derivative of the logarithm gives
$$\frac{G'_d(s)}{G_d(s)}-\frac{G'_d(2s)}{G_d(2s)}=\frac{\zeta'(2s)}{\zeta(2s)}+\frac{T_d'(2s)}{T_d(2s)}-\frac{L'(s,\chi_d)}{L(s,\chi_d)}.$$
Thus
\begin{align}\frac{G'_d(s)}{G_d(s)}-\frac{G'_d(2^ms)}{G_d(2^ms)}&=
\sum_{k=0}^{m-1}\left (
\frac{G'_d(2^ks)}{G_d(2^ks)}-\frac{G'_d(2^{k+1}s)}{G_d(2^{k+1}s)}\right
)\\
&=\sum_{k=1}^{m}\frac{\zeta'(2^ks)}{\zeta(2^ks)}+\sum_{k=1}^{m}\frac{T_d'(2^ks)}{T_d(2^ks)}-
\sum_{k=0}^{m-1}\frac{L'(2^ks,\chi_d)}{L(2^ks,\chi_d)}.
\end{align}
 By \eqref{G} we have that for
$s\ge  2$, $|\frac{\zeta'(s)}{\zeta(s)}|\le \sum_{n\ge
2}\frac{\Lambda(n)}{n^{s}-1}\le \frac{\log 2}{2^s-1}+\sum_{n\ge
3}\frac{\log n}{n^{s}-1}\le  \frac 43\frac{\log 2}{2^s}+\frac
98\sum_{n\ge 3}\frac{\log n}{n^{s}}\le \frac 43\frac{\log
2}{2^s}+\frac 98\int_2^{\infty}\frac{\log x}{x^{s}}dx=\frac
43\frac{\log 2}{2^s}+\frac 98\Bigl (\frac{\log
2}{2^{s-1}(s-1)}+\frac 1{2^{s-1}(s-1)^2}\Bigl )\le \frac
1{2^s(s-1)}\Bigl (\frac{20\log 2+8}9\Bigr )\le \frac{5}2\cdot
\frac{2^{-s}}{s-1}.$ Thus, $|\frac{\zeta'(2^k)}{\zeta(2^k)}|\le
\frac{5}2\cdot \frac{2^{-2^k}}{2^k-1}.$ The same estimate holds for
$|\frac{G_d'(2^k)}{G_d(2^k)}|$, $|\frac{T_d'(2^k)}{T_d(2^k)}|$ and
$|\frac{L'(2^k,\chi_d)}{L(2^k,\chi_d)}|$. When $m\to \infty$ and
then $s\to 1$ we get
\begin{align}
\sum_{p}\frac{(d/p)\log
p}{p-1}=\sum_{k=1}^{\infty}\frac{\zeta'(2^k)}{\zeta(2^k)}
-\sum_{k=0}^{\infty}\frac{L'(2^k,\chi_d)}{L(2^k,\chi_d)}+\sum_{k=1}^{\infty}\frac{T'_d(2^k)}{T_d(2^k)}.
\end{align}
Finally we observe that $\frac{T'_d(2^k)}{T_d(2^k)}=\sum_{p\mid
d}\frac{\log p}{p^{2^k}-1}$, so
$\sum_{k=1}^{\infty}\frac{T'_d(2^k)}{T_d(2^k)}=\sum_{p\mid d}s_p$.
\end{proof}

The advantage of the lemma above is that the series involved
converge very fast. For example,
$\sum_{k=0}^{\infty}\frac{L'(2^k,\chi_d)}{L(2^k,\chi_d)}=\sum_{k=0}^{6}\frac{L'(2^k,\chi_d)}{L(2^k,\chi_d)}+Error$
with $|Error|\le 10^{-40}$.

Hence we can write $B_f=C_0+C_d+C(f)$ where $C_0$ is an universal
constant, $C_d$ depends only on $d$, and $C(f)$ depends on $f$. More
precisely,
\begin{align*}
C_0&=\gamma -1-2\log
2-\sum_{k=1}^{\infty}\frac{\zeta'(2^k)}{\zeta(2^k)}=-1.1725471674190148508587521528364 \\
C_d&=\sum_{k=0}^{\infty}\frac{L'(2^k,\chi_d)}{L(2^k,\chi_d)}-\sum_{p\mid
d}s_p\\
C(f)&=\frac 1{\phi(q)}\sum_{\substack{1\le r\le
q\\(r,q)=1}}\log\Bigl (1+\frac rq\Bigr )+\log a+\sum_{p\mid 2aD}\log
p\Bigl(\frac{1+(d/p)}{p-1}- \sum_{k\ge 1}\frac{s(f,p^k)}{p^k}\Bigr
).
\end{align*}

The values of $s_p$ and $\sum_{k\ge 0}L'(2^k,\chi_d)/L(2^k,\chi_d)$,
can be calculated with MAGMA with high precision. We include some
values of $s_p$,
\begin{align*}
s_2&=0.279987673370859807200459206376\dots\\
s_3&=0.151226686598727076356318275233\dots\\
s_5&=0.069643260624011195267442944307\dots\\
s_7&=0.041350928217815118656218939260\dots\\
s_{17}&=0.009871469313243775687197132626\dots\\
\end{align*}
some values of $C_d$,
\begin{align*}{}&C_{-4}&=&+0.346538435736895987549
-s_2&=&+0.066550762366036180349\\
&C_{-8}&=&-0.076694093066485311184-s_2&=&
-0.356681766437345118384\\
&C_{8}&=&+0.809903104673738787384-s_2&=&+0.529915431302878980184\\
&C_{-3}&=&+0.586272400297149523649-s_3&=&+0.435045713698422447292\\
&C_5&=&+
1.172449603551261794528-s_5&=&+1.102806342927250599260\\
&C_{-7}&=&-0.070022837990444988815
-s_7&=&-0.111373766208260107471\\
&C_{12}&=&+0.564588639325865961984-s_2-s_3
&=&+0.133374279356279078427\\
&C_{-15}&=&-0.486320692903261758405-s_3-s_5
&=&-0.707190640126000030028\\
&C_{17}&=&+0.289109343784025529610-s_{17}&=&+0.279237874470781753922
\end{align*}
and some values of $C(f(x))$:
\begin{align*}
&C(x^2+1)&=&\ (3\log 2)/2&=&\ 1.039720770839917964125 \dots\\
&C(x^2+2)&=&\ (3\log 2)/2&=&\ 1.039720770839917964125\dots\\\
&C(x^2-2)&=&\ (3\log 2)/2&=&\ 1.039720770839917964125\dots\\\
&C(x^2+x+1)&=&\ \log 2+(\log 3)/6&=&\ 0.876249228671296924649\dots\\\
&C(x^2+x-1)&=&\ \log 2+(\log
5)/(20)&=&\ 0.773619076181650328147\dots\\\
&C(x^2+x+2)&=&\ \log 2+(\log 7)/(42)&=&\
0.739478374585071816681\dots\\\
&C(x^2+2x-2)&=&\ (3\log 2)/2+(\log
3)/6&=&\ 1.222822818951269579358\dots\\\
&C(2x^2+1)&=&\ 3\log 2&=&\ 2.079441541679835928251\dots\\\
&C(2x^2-1)&=&\ 3\log 2&=&\ 2.079441541679835928251\dots\\\
&C(2x^2+x+1)&=&\ 2\log 2+\log 3+(\log
7)/(42)&=&\ 1.838090663253181508076\dots\\\
&C(2x^2+x+2)&=&\ \log 2+(7\log 3)6+(\log 5)/(20)&=&\
2.055333412961111634775\dots\\\
&C(2x^2+x-2)&=&\ \log 2+\log 3+(\log(17))/(272)&=&\
1.802175694757673442283\dots\\\
&C(2x^2-x+1)&=&\ \log 2+\log 3+(\log(7))/(42)&=&\
1.838090663253181508076
\\\
&C(2x^2-x+2)&=&\ \log 2+(7\log 3)/6+(\log 5)/(20)&=&\
2.055333412961111634775\dots
\\\
&C(2x^2-x-2)&=&\ \log 2+\log 3+(\log(17))/(272)&=&\
1.802175694757673442283\dots
\\\
&C(2x^2+2x+1)&=&\ 3\log 2&=&\ 2.079441541679835928251\dots
\\\
&C(2x^2+2x-1)&=&\ 3\log 2+(\log 3)/6&=&\
2.262543589791187543484\dots
\end{align*}
Table below contains the constant $B=B_f$  for all irreducible
quadratic polynomial $f(x)=ax^2+bx+c$ with $0\le a,|b|,|c|\le 2$.
When $f_1,f_2$ are irreducible quadratic polynomials such that
$f_1(x)=f_2(x+k)$ for some $k$, we only include one of them since
$L_n(f_1)=L_n(f_2)+O(\log n)$.
\begin{center}
\begin{tabular}{|c|c|c|l|}  \hline
$f(x)$ & $d$ & $q$ & $\qquad \qquad \qquad B_f$
\\\hline $x^2+1$ & -4 & 1 &\ - 0.06627563421306070638\dots
\\\hline $x^2+2$ & -8& 1 &\  - 0.48950816301644200511\dots
\\\hline $x^2-2$& 8& 1 &     + 0.39709034723782093451\dots
\\\hline $x^2+x+1$& -3& 1&   + 0.13874777495070452108\dots
\\\hline $x^2+x-1$& 5 & 1&   + 0.70387825168988607654\dots
\\\hline $x^2+x+2$&-7& 1&\   - 0.54444255904220314164\dots
\\\hline $x^2+2x-2$& 8& 1&   + 0.18364993088853380692\dots
\\\hline $2x^2+1$&-8&1&      + 0.55021260782347595900\dots
\\\hline $2x^2-1$&8& 1&      + 1.43680980556370005757\dots
\\\hline $2x^2+x+1$&-7&2&    + 0.55416972962590654974\dots
\\\hline $2x^2+x+2$&-15&2&   + 0.17559560541609675388\dots
\\\hline $2x^2+x-2$&17& 2&   + 0.90886640180944034534\dots
\\\hline $2x^2-x+1$&-7& 2&   + 0.55416972962590654974\dots
\\\hline $2x^2-x+2$&-15 &2 & + 0.17559560541609675388\dots
\\\hline $2x^2-x-2$ &17 &2 & + 0.90886640180944034534\dots
\\\hline $2x^2+2x+1$ &-4 &1 &+ 0.97344513662685725774\dots
\\\hline $2x^2+2x-1$ &12 &1 &+ 1.22337070172845177105\dots
\\\hline
\end{tabular}
\end{center}

\begin{center}

\medskip

 Table below
shows the error term $E_f(n)=\log L_n(f)-n\log n-B_fn$ for the
polynomials above and some values of $n$.

\smallskip

\begin{tabular}{|c|c|c|c|c|
c|c|c|c|}
  \hline
 $f(x)$ & $x^2+1$ & $x^2+2$ &$x^2-2$ &$x^2+x+1$ &$x^2+x-1$ &$x^2+x+2$ & $x^2+2x-2$\\
  \hline
  $E_f(10^2)$ & -18 &-36  &-7 &-6 & -12   &+9 &-17 \\ \hline
  $E_f(10^3)$ & +6 &-11  & -46& -9 & -91  &-20 &-97 \\ \hline
  $E_f(10^4)$ & -111 &-263  &-54 &+17 &-208 &-218 &-297 \\ \hline
  $E_f(10^5)$ & +34 &-761  & -466&-654  &-253 &-2120 &-553 \\ \hline
  $E_f(10^6)$ & -2634 &-1462 &-764 &-2528  &-1075 & +687&-454 \\ \hline
 $E_f(10^7)$ & -1557 &-8457 &-1472 &-1685 &-9636 &-686 &-6336 \\ \hline
\end{tabular}

\

\begin{tabular}{|c|c|c|c|c|
c|c|c|c|}
  \hline
 $f(x)$ & $2x^2+1$ & $2x^2-1$ &$2x^2+x+1$ &$2x^2+x+2$ &$2x^2+x-2$ &$2x^2-x+1$ \\
  \hline
  $E_f(10^2)$ &-15  &-19 &-1 &-34 &-5 &-22  \\ \hline
  $E_f(10^3)$ & -1 & -69& +6& +4 &-37 & -126 \\ \hline
  $E_f(10^4)$ &-301  &-233 & +18& -295&-198 &-43  \\ \hline
  $E_f(10^5)$ & -251 &-182 &-1289 & +27&-1193 &+177  \\ \hline
  $E_f(10^6)$ & +1084 &-159 & +235& +1169&-4856 &-3077  \\ \hline
 $E_f(10^7)$ & -14821 & -10525&-2553 & +1958&-16758 &-5459  \\ \hline
\end{tabular}

\

\begin{tabular}{|c|c|c|c|c|c|
c|c|c|c|}
  \hline
 $f(x)$ &$2x^2-x+2$& $2x^2-x-2$ & $2x^2+2x+1$ &$2x^2+2x-1$  \\
  \hline
  $E_f(10^2)$ &-5  &-17 &-9 & -14\\ \hline
  $E_f(10^3)$ &-123  &-18 &-89 &  -41\\ \hline
  $E_f(10^4)$ &+74  &-136 & +9&  +58\\ \hline
  $E_f(10^5)$ &-2083  &-516 &-232 & -331 \\ \hline
  $E_f(10^6)$ &-4851  &+3532 &-2876 & -931 \\ \hline
 $E_f(10^7)$ &-18152  &+907 &-10624 & +689 \\ \hline
\end{tabular}
\end{center}

\section{Quadratic reducible polynomials} To complete the problem of
estimating the least common multiple of quadratic polynomials we
will study here the case of reducible quadratic polynomials. Being
this case much easier than the irreducible case, we will give a
complete description for sake of the completeness.

If $f(x)=ax^2+bx+c$ with $g=(a,b,c)>1$, it is easy to check that
$\log L_n(f)=\log L_n(f')+O(1)$ where $f'(x)=a'x^2+b'x+c'$ with
$a'=a/g,\ b'=b/g,\ c'=c/g$.

If $f(x)=(ax+b)^2$ with $ (a,b)=1$ then, since $(m^2,n^2)=(m,n)^2$,
we have that $L_n((ax+b)^2)=L_n^2(ax+b)$ and we can apply \eqref{Ba}
to get
\begin{equation}
\log \text{l.c.m.}\{(a+b)^2,\dots ,
(an+b)^2\}\sim 2n\frac{a}{\phi (a)}\sum_{\substack{1\le k\le a\\
(k,a)=1}}\frac 1k.
\end{equation}

Now we consider the more general case $f(x)=(ax+b)(cx+d),\
(a,b)=(c,d)=1$.
\begin{thm} Let $f(x)=(ax+b)(cx+d)$ with $(a,b)=(c,d)=1$ and $ad\ne bd$. Let $q=ac/(a,c)$. We have
\begin{equation}
\log L_n(f)\sim \frac{n}{\varphi(q)}\sum_{\substack{1\le r\le q,\
(r,q)=1}}\max \Bigl (\frac a{(br)_a},\frac c{(dr)_c}\Bigr
).\end{equation}
\end{thm}
\begin{proof}
Suppose $p^2\mid L_n(f)$. It implies that $p^2\mid (ai+b)(ci+d)$ for
some $i$. If $p\mid ai+b$ and $p\mid ci+d$ then $p\mid (ad-bc)i$. If
$p\nmid (ad-bc)$ then $p\mid i$ and consequently $p\mid b$ and
$p\mid d$. Thus, if $p\nmid (ad-bc)bd$ and $p^2\mid (ai+b)(ci+d)$
then $p^2\mid (ai+b)$ or $p^2\mid (ci+d)$. In these cases $p\le
M_n=\max (\sqrt{an+b},\sqrt{cn+d},|(ad-bd)bd|)$.

Thus we write
\begin{equation}
L_n(f)=\prod_{p\le
M_n}p^{\beta_p(n)}\prod_{p>M_n}p^{\epsilon_p(n)}=\prod_{p\le
M_n}p^{\beta_p(n)-\epsilon_p(n)}\prod_{p}p^{\epsilon_p(n)},
\end{equation}
where $\epsilon_p(n)=1$ if $p\mid f(i)$ for some $i\le n$ and
$\epsilon_p(n)=0$ otherwise.
 Since $p^{\beta_p(n)}\le f(n)$ we have that
$\beta_p(n)\ll \log n/\log p$ and then
\begin{equation}
\sum_{p\le M_n}(\beta_p(n)-\epsilon_p(n))\log p\ll (\log
n)\pi(M_n)\ll \sqrt n.
\end{equation}
Thus,
\begin{equation}
\log L_n(f)=\sum_{\substack{p\mid f(i)\\\text{for some } i\le
n}}\log p+O(\sqrt n)
\end{equation}
Let $q=ac/(a,c)$. Suppose that $p\equiv r^{-1}\pmod{q},\ (r,q)=1$.
Let $k=(br)_a$ the least positive integer such that $k\equiv br\pmod
a$. Then $p\mid (ai+b)$ for some $i\le n$ if and only if $kp\le
an+b$. Similarly, let $j=(dr)_c$ be the least positive integer such
that $j\equiv dr\pmod c$. Again, $p\mid (ci+d)$ for some $\le i\le
n$ if $jp\le cn+d$. Thus, the primes $p\equiv r^{-1}\pmod{ac}$
counted in the sum above are those such that $p\le
\max(\frac{an+b}{k},\frac{cn+d}{j})$. The prime number theorem for
arithmetic progressions implies that there are $\sim \frac
n{\varphi(q)}\max(\frac a{k},\frac c{j})$ of such primes.

We finish the proof summing up in all $1\le r\le q,\ (r,q)=1$.
\end{proof}

\section{Some remarks about the error term}

It is known that the estimate $E(n)=\log \text{l.c.m.}\{1,\dots,
n\}-n=O(n^{1/2+\epsilon})$ is equivalent to the Riemann hypothesis.
Probably it is also true that $E_f(n)=O(n^{1/2+\epsilon})$ for any
irreducible quadratic sequence, but  it is clear that to prove that
$E_f(n)=O(n^{\theta})$ for some $\theta<1$ is a very hard problem.

Recently K. Homma \cite{K} has proved that if $D<0$ then
$$\# \{\nu/p\in I:\ 0<\nu<p\le x:\ f(\nu)\equiv 0\pmod
p\}=\pi(x)\Bigl (1+O(1/(\log x)^{\theta})\Bigr)$$ for any
$\theta<8/9$. Using this result and the known error term for the
prime number theorem for arithmetic progressions it is possible to
prove that $E_f(n)=O(\frac n{\log^{\alpha}n})$ for some $\alpha>0$,
when $f(x)$ is an irreducible quadratic polynomial of the form
$f(x)=ax^2+c,\ a,c>0$.


%
%
%
%
%

%
%
%

\end{document}